---------------------------------------------------------------------------
############################## .TEX FILE ###################################
---------------------------------------------------------------------------

\documentstyle[aps]{revtex}            
\begin{document}

\thispagestyle{empty}


\draft
\preprint{}
\title{Local Hypercomplex Analyticity}
\author{Stefano De Leo$^{\; (a,b,c)}$
        and
        Pietro Rotelli$^{\; (a)}$ }
\address{$^{(a)}$ Dipartimento di Fisica, 
                Universit\`a degli Studi Lecce and INFN, 
                Sezione di Lecce\\
                via Arnesano, CP 193, 73100 Lecce, Italia\\
          $^{(b)}$ Dipartimento di Metodi e Modelli Matematici per le Scienze
                Applicate,\\
                via Belzoni 7, 35131 Padova, Italia\\ 
          $^{(c)}$ Instituto de Matem\'atica, Estat\'{\i}stica e 
                   Computa\c{c}\~ao Cient\'{\i}fica, UNICAMP,\\
                Cidade Universit\'aria ``Zeferino Vaz'', 
                CP 6065, 13081-970, Campinas, S.P., Brasil} 
\date{February, 1997}
\maketitle


\begin{abstract}

The notion of analyticity is studied in the context of hypercomplex 
numbers. A critical review of the problems arising from the conventional 
approach is given. We describe a local analyticity condition which 
yields the desired type of hypercomplex solutions. The result is the 
definition of a generalized complex analyticity to hypercomplex space.


\end{abstract}


\section{Introduction}
\label{s1}

When in 1843, searching for a new numerical field  that in three dimensions 
operates in 
analogy to complex numbers in two dimensions, the Irish mathematician 
William Rowan Hamilton discovered quaternions~\cite{H43}, he probably was not 
aware of what difficulties hide behind quaternions and what hostilities would
arise in the Science community against such non commutative numbers.

With the notable exception of Hamilton himself, quaternions have always met 
difficulties in finding a   
physical incarnation~\cite{G83}. 
Notwithstanding their (implicit) technical use in the 
{\sl Treatise of Electricity and Magnetism} by Maxwell~\cite{M},  where
the $\nabla$-operator is given in terms of the three quaternionic imaginary 
units $i$, $j$, $k$
\[ \nabla = i \partial_x + j \partial_y + k \partial_z~,\]
quaternions were
quickly moved aside by the subsequent and practical vectorial algebra
formulated of Gibbs and Heaviside~\cite{NAH} in the 1880's. 
Gibbs and Heaviside,  
starting from quaternionic products, introduced the usual three dimensional 
scalar and vector product. It was these two products which were 
then preferentially applied in 
Physics, in alternative to the quaternionic product. Thus, quaternions fell in 
disuse in Physics not many decades after they were introduced. 

While the intrinsic non commutativity of Quantum Mechanics should have 
encouraged the use of the quaternionic field as the natural underlying 
numerical 
field~\cite{FW}, the difficulties encountered in 
quaternionic calculations, together with the natural 
conservativism of the Science community have discriminated against the 
application of Hamilton's number field. In particular, problems arise in the 
appropriate definitions of tensor product spaces, in  
the corrected definition of determinant and transpose for quaternionic 
matrices, in the reproduction of standard trace theorems, and in the passage 
from hermitian to antihermitian operators, etc. 

Nevertheless, after the fundamental paper of Horwitz and Biedenharn on Gauge 
Theories
\& Quantization~\cite{HB84} and Adler's work on Quaternionic
Quantum Mechanics \& Quantum Fields~\cite{AB,AP}, 
quaternions have gained new impulse  in Quantum Physics.  Much recent 
progress in quaternionic applications to Mathematics and Physics have been 
obtained  and thus  
quaternions are considered with more interest nowadays. 
Among the recent applications of  
quaternions we recall quaternionic formulations for 
Tensor Products~\cite{TP}, Group Representations~\cite{GR}, Special
Relativity~\cite{SR}, Dark Matter~\cite{DM}, Relativistic 
Quantum Mechanics~\cite{R89,RE}, 
Scattering \& Decay Problems~\cite{SD}, 
Field Theory~\cite{FT}, 
Hidden Symmetry~\cite{HS}, Electroweak Model~\cite{EM}, 
GUTs~\cite{GUT}, Preonic Model~\cite{PM}, etc. Experimental tests, which
might highlight the possible quaternionic nature of wave functions, have 
also been proposed~\cite{ET}.

Within Mathematics itself an important question remains pertinent: 
How is the well known complex analyticity condition,
\[
\partial_{x_0} f(z,\bar{z}) = - i \partial_{x_1} f(z,\bar{z}) \quad \quad \quad
[~z= x_0 + ix_1 ~~~(x_0, x_1 \in {\cal R})~]~,
\]
best  extended to hypercomplex numbers including quaternions?

In this paper, we start with  a critical review of the conventional approach 
to quaternionic analysis, showing the disadvantages of defining a 
quaternionic derivative in the standard form
\[
\partial_q \equiv \alpha \partial_{x_0} + \beta \partial_{x_1} + \gamma
\partial_{x_2} + \delta \partial_{x_3}~,
\]
with $\alpha, \; \beta, \gamma, \; \delta $ ``constant'' quaternionic numbers.

The main goal of our work is the demonstration that a new   
generalization of complex analyticity to quaternionic numbers can be 
easily obtained. This occurs thanks to a local representation 
of the complex imaginary unit, that yields the possibility of defining a  
{\em local quaternionic derivative} which overcomes some previous problems. 
Such a generalization is easily extended to octonionic numbers.

This paper is structured as follows: In Section II, we introduce  
the quaternionic algebra and  the {\em barred operators}. 
In Section III, after a brief introduction to 
complex differentiability we discuss the standard quaternionic approaches. 
In Section IV, we recall some of the 
essential aspects of Fueter's calculus.   
Section V is directed to the study of ``global''  
analyticity conditions. In this Section we show the 
impossibility of finding a suitable  quaternionic derivative by adopting 
the standard viewpoint.
Section VI deals with the ``local'' hypercomplex version 
of the Cauchy-Riemann equations. Our conclusions are drawn in the final 
Section.


\section{Quaternionic Algebra and Barred Operators}
\label{s2}

Let us introduce the quaternionic algebra by trying to follow the conceptual 
approach of Hamilton. We begin by looking for numbers of the form 
$x_1+ix_2+jx_3$, 
with $i^2=j^2=-1$, which will do for three-dimensional space what complex 
numbers have done 
for the plane. Influenced by the existence of a complex number norm
\[ \bar{z} z = (\mbox{Re}z)^2 +  (\mbox{Im}z)^2~,\]
when we look at its generalization 
\[ (x_1-ix_2-jx_3) (x_1+ix_2+jx_3) = 
    x_1^2+x_2^2+x_3^2 - (ij+ji) x_2 x_3~,\]
to obtain a real number, it is necessary to adopt  the anticommutative law of 
multiplication for the 
imaginary units. We now note that 
with only two imaginary units we have no chance of constructing a new 
numerical field, because assuming
\begin{eqnarray*}
ij=\alpha_1 + i \alpha_2 + j \alpha_3~,\\
ji=\beta_1 + i \beta_2 + j \beta_3~,\\
\end{eqnarray*}
and
\[ ij=-ji~,\]
we find
\[ \alpha_{1,2,3}= \beta_{1,2,3}=0~.\]
Thus we must introduce a third imaginary unit $k\neq i\neq j$, with
\[ k=ij~. \] 
This noncommutative field is therefore characterized by three imaginary 
units $i$, $j$, $k$ which satisfy the following multiplication rules
\begin{equation}
i^2=j^2=k^2=ijk=-1~.
\end{equation}
Numbers of the form
\begin{equation}
q=x_0+ ix_1+jx_2+kx_3~~~~~~~~(~x_{0,1,2,3} \in {\cal R}~)~,
\end{equation}
are called (real) {\em quaternions}. They 
are added, subtracted and multiplied according to the usual laws of 
arithmetic, except for the commutative law of multiplication.

An important difference between quaternionic and complex numbers is related 
to the definition of the conjugation operation. Whereas with complex 
numbers we can define only one type of conjugation
\[ i \rightarrow -i~,\]
working with quaternionic numbers we can introduce 
different conjugation operations.  
Indeed, with three imaginary units we
have the possibility to define besides the standard conjugation
\begin{equation}
\bar{q} = x_0 - ix_1 - jx_2 -kx_3~,
\end{equation}
 the six new operations
\begin{eqnarray*}
(i,j,k) & ~\rightarrow ~& (-i,+j,+k)~,~(+i,-j,+k)~,~(+i,+j,-k)~;\\
(i,j,k) & ~\rightarrow ~& (+i,-j,-k)~,~(-i,+j,-k)~,~(-i,-j,+k)~.
\end{eqnarray*}
These last six conjugations can be concisely expressed in terms of 
$q$ and $\bar{q}$ by
\begin{eqnarray*}
q & ~\rightarrow ~& -i\bar{q}i~,~-j\bar{q}j~,~-k\bar{q}k~,\\
q & ~\rightarrow ~& -iqi~,~-jqj~,~-kqk~.
\end{eqnarray*}
It might seem that the only independent conjugation be represented by 
$\bar{q}$. Nevertheless, $\bar{q}$ can also be expressed 
in terms of  $q$ 
\begin{equation}
\label{qbar}
\bar{q} = -\frac{1}{2} \; (q+iqi+jqj+kqk)~.
\end{equation}
So we conclude that while within a complex theory $z$ and $\bar{z}$ represent 
independent variables, in a quaternionic theory all the seven 
conjugations can be defined as involutions of the quaternionic variable $q$. 
We will return to this point later. 

Let us now introduce formally  the so-called {\em barred operators}, 
already implicit in the previous formulas.  
In recent years the left/right-action of the quaternionic imaginary
units, expressed by barred operators~\cite{DRT}, has been very
useful in overcoming difficulties owing to the noncommutativity of
quaternionic numbers. Among the successful applications of barred
operators we mention the one-dimensional quaternionic formulation of
Lorentz boosts~\cite{SR} and new possibilities for Quaternionic Group 
Theories~\cite{GUT}. Partially-barred quaternions also appear in
Quantum Mechanics, e.g they allow an appropriate definition of
the momentum operator, in a quaternionic version of Dirac's
equation~\cite{R89} and Electroweak Model~\cite{EM}.

The most general barred quaternionic operator is represented by
\begin{equation}
\label{bo}
q_0 + q_1 \mid i + q_2 \mid j + q_3 \mid k \quad \quad \quad
(~q_{0,1,2,3} \in {\cal H}~)~.
\end{equation}
where $1\mid i$, $1\mid j$ and $1\mid k$ indicate the right-action of
the quaternionic imaginary units. The application of $q_1\mid q_2$ to a 
third quaternion $q_3$ is by definition $q_1q_3q_2$. 
The multiplication rules are simple
once seen. For example, in the multiplication of two particular 
barred quaternions, $j\mid k$ and $i\mid j$, we find
\[ (j\mid k) \times (i\mid j) = ji \mid jk = - k\mid i~.\]
We also observe that in this example $j\mid k$ and $i\mid j$ represent 
commuting quaternionic barred operators.

Let us now analyze the link between complex/quaternion numbers and real 
matrices. We know that the action of a complex number on $z$ can be 
given in terms of real matrices. In fact, identifying $z$ with the 
vector column
\begin{equation}
\label{vec}
\left( \begin{array}{c} x_0\\ x_1 \end{array} \right)~,
\end{equation}
the action on $z$ by the complex number $a+ib$ is expressed by the 
following real matrix
\begin{equation}
\label{mc} 
\left( \begin{array}{cc} a & -b \\ b & a \end{array} \right)~. 
\end{equation}
Obviously this involves only particular real matrices. 
This corresponds to the fact that complex numbers are characterized 
by {\em two} real parameters whereas the most general $2\times 2$ real 
matrix by four. The restriction in Eq.~(\ref{mc}) 
precludes one from obtaining from~(\ref{vec}) 
the complex conjugate
\[ \bar{z} \equiv \left( \begin{array}{c} x_0\\ $-$x_1 \end{array} \right)~, \]
and whence we see again, in another way, that $z$ and $\bar{z}$ represent 
independent variables. The situation, drastically changes for quaternions. 
If we represent a 
quaternionic number, $q$, by the column vector
\[ \left( \begin{array}{c} x_0\\ x_1\\ x_2 \\ x_3 \end{array} \right)~, \]
the most general action we can perform on this column vector is obviously 
given by a $4\times 4$ real matrix ($16$ real parameters). But, due to the 
noncommutative nature 
of quaternions, the most general transformation on quaternionic numbers is 
represented by the barred quaternions~(\ref{bo}), characterized by 
just $16$ real parameters. So we have a {\em full} connection between 
$4\times 4$ real matrices and quaternions. 
The situation can be concisely summarized as 
follows
\begin{center}
\begin{tabular}{|cccc|}\hline
~~~~Field~~~~    & 
~~~~State's ~~~~ &
~~~~Operator's ~~~~& 
~~~~Real Matrices~~~~\\
   & 
~~~~~Real Parameters~~~~~ &
~~~~~Real Parameters~~~~~ & \\ \hline
$\cal C$ &        two   &    two                 & $O(2) \in Gl(2)$\\
$\cal H$ &        four   &   sixteen                & $Gl(4)$\\ \hline
\end{tabular}
\end{center}

There is another important distinction between complex and quaternionic 
numbers. The requirement that a function of a complex variable 
$z=x_0+ix_1$ should 
be a complex polynomial function $f(x_0,x_1)+ig(x_0,x_1)$ picks out a 
proper subset of the polynomial functions. However, the 
corresponding requirement of a function of a quaternionic variable, $q$, 
namely that it should be a sum of monomials of the type  
\[a_oqa_1...a_{z-1}qa_r~~~~~~~~(a_n \in {\cal H})~.\] 
places {\em no} restriction 
on the function. In contrast to the complex case the 
coordinates $x_0,\vec{x}$ can themselves be written as quaternionic 
polynomials
\begin{eqnarray*}
x_0  &~~=~~& P_r \, q\equiv \frac{1-i\mid i - j \mid j - k \mid k}{4} \, q~,\\
(~ix_1, \, jx_2, \, kx_3~) &~~=~~& 
(~-iP_ri \, q, \,-jP_rj \, q, \, -kP_rk \, q~) \equiv 
(~P_i \, q~, \, P_j \, q~, \, P_k \, q~)~,
\end{eqnarray*}
or explicitly,
\begin{eqnarray*}
x_0  &~~=~~& ~~\frac{1}{4} ~  (1-i\mid i - j \mid j - k \mid k) \, q~,\\
x_1  &~~=~~& - \frac{i}{4} ~  (1-i\mid i + j \mid j + k \mid k) \, q~,\\
x_2  &~~=~~& - \frac{j}{4} ~  (1+i\mid i - j \mid j + k \mid k) \, q~,\\
x_3  &~~=~~& - \frac{k}{4} ~  (1+i\mid i + j \mid j - k \mid k) \, q~,
\end{eqnarray*}
and so every real polynomial in $(x_0, \vec{x})$ can be written as 
quaternionic polynomial in $q$ by combinations of
\[ q~~,~~-iqi~~,~~-jqj~~,~~-kqk~~.\]
Thus, a theory of quaternionic power series will be the same as a 
theory of real power series  functions in ${\cal R}^4$. 

Consequently, while working 
within complex numbers we have no possibility to express $\bar{z}$ by $z$, 
so that it is necessary to consider a general function $f$ as  
$f(z,\bar{z})$,
\[ f_{0}(x_0,x_1) + i f_{1}(x_0,x_1) \rightarrow f(z,\bar{z})
~~~~~~~~(~f_{0,1} \in {\cal R}~,~f \in {\cal C}~)~.\] 
In quaternionic theories, due to the noncommutativity,   
$\bar{q}$ (and all other conjugates) can be written in terms of $q$, 
by quaternionic barred operators, e.g. 
\[ \bar{q} = - \frac{1+i\mid i + j \mid j + k \mid k}{2} \, q~,\] 
thus, we can in all generality  refer to $f(q)$ instead of 
$f(q,\bar{q},...)$,
\[ f_{0}(x_0,\vec{x}) + i f_{1}(x_0,\vec{x}) +jf_{2}(x_0,\vec{x}) 
+ k f_{3}(x_0,\vec{x}) 
\rightarrow f(q)
~~~~~~~~(~f_{0,1,2,3} \in {\cal R}~,~f \in {\cal H}~)~.\]


\section{Quaternionic Differentiability}
\label{s3}

It is well known that a point in the two-dimensional real Euclidean space
$(x_0, x_1)$ can be naturally expressed by a complex number
\[ z=x_0 + i x_1 \]
and that the most general conformal coordinate transformation in two
dimension is analytic in this complex coordinate. In four dimensions, 
coordinates can be represented by a quaternion as
naturally as, in two dimensions, they are unified into a complex
number. Thus, it is natural to ask whether the notion of
analyticity can be extended to quaternions.

Before discussing hypercomplex analyticity, let us briefly recall the
traditional approach to the notion of complex analyticity. A function whose 
arguments 
are complex numbers $x_0 \pm i x_1$ is defined 
as
\begin{equation}
\label{complex}
f(z,  \bar{z})= u(x_0,  x_1) + i v(x_0,  x_1)~,
\end{equation}
where $u(x_0,  x_1)$ and $v(x_0,  x_1)$ are real functions of the real 
arguments $x_0$ and $x_1$. The standard theory of analytic functions deals 
only with a restricted but important class of functions of two real 
variables, namely, those functions that satisfy the requirement 
of being ``differentiable''. Though differentiability of a function of a 
complex variable places, as we said,  a severe 
limitation on the functions one is allowed to consider, it leads 
to a theory of these functions that is both elegant and extremely 
powerful. We aim to extend the standard discussion on the differentiability 
of functions of a complex variable to quaternionic functions of 
a {\em quaternionic} (and even octonionic)  
variable.

The derivative of a function of a complex variable with respect to the 
{\em argument} $z$ is formally defined in the same way as it is for 
functions of a real variable (complex numbers commute and so the position of 
the factor $1/\Delta z$ below is not relevant)
\begin{equation}
\label{derz}
\partial_z  f(z,  \bar{z}) = 
\lim_{\Delta z \rightarrow 0}
\frac{f(z+\Delta z, \bar{z}+ \Delta \bar{z}) - f(z)}{\Delta z}~,
\end{equation}
with $\Delta z=\Delta x_0 + i \Delta x_1$. 
The limit in the previous equation will, in general, depend upon the order on 
which $\Delta x_0$, $\Delta x_1$ tend to zero. We say that a function of a 
complex argument 
is {\em differentiable} if the limit~(\ref{derz}) exists, is finite, and 
does not depend on the manner in which one takes the limit. Whereas in the 
case of a real variable, one can approach a given point only in two ways 
(either from the left or from the right), a point in the complex plane can be 
reached from an infinite number of directions. Thus, instead of the single  
requirement for a function of a real variable,  
\[ Left~~\mbox{limit} ~~~~\equiv ~~~~Right~~\mbox{limit}~, \]
 an infinite number of such 
requirements has to be satisfied in order to ensure the differentiability 
of a function of a complex variable. One can understand therefore, why the 
property of being differentiable is, in the case of functions of a complex 
variable, very much more restrictive than it is for functions of a real 
variable.

The {\em necessary} and {\em sufficient} condition for a function 
of a complex variable to be 
differentiable at a 
given point is that, at that point, it obeys the Cauchy-Riemann conditions. 
Such conditions are obtained by imposing the requirement that the 
right-hand side of~(\ref{derz}) yields the same result, whatever the order 
in which the limit $\Delta x_0 \, \Delta x_1 \rightarrow 0$ is taken. By 
first setting $\Delta x_1=0$ and then taking the limit 
$\Delta x_0 \rightarrow 0$ we find
\begin{equation}
\label{1a}
\partial_z  f(z,  \bar{z}) =
\partial_{x_0} u(x_0,  x_1) + i \partial_{x_0} v(x_0, x_1)~.
\end{equation}
In the other case, in which we first set $\Delta x_0=0$ and then take the 
limit $\Delta x_1 \rightarrow 0$ we find
\begin{equation}
\label{1b}
\partial_z  f(z,  \bar{z}) =
-i \partial_{x_1} u(x_0,  x_1) +  \partial_{x_1} v(x_0, x_1)~.
\end{equation}
By equating the real and imaginary parts of Eqs.~(\ref{1a}, \ref{1b}), 
we obtain the Cauchy-Riemann conditions:
\[ \partial_{x_0} u = \partial_{x_1} v
~~~~~\mbox{and} ~~~~~ 
\partial_{x_1} u  = -\partial_{x_0} v~,
\]
or equivalently
\begin{equation}
\label{zzx}
\partial_{x_0} f(z,\bar{z}) = - i \partial_{x_1} f(z,\bar{z})~.
\end{equation}
The real and imaginary parts of a differentiable function separately 
satisfy the Laplace equation
\[ ( \partial_{x_0}^2 + \partial_{x_1}^2)     u =
( \partial_{x_0}^2 + \partial_{x_1}^2)     v=0~,
\]
and are therefore {\em harmonic} functions of two variables (the converse 
is {\em not} necessarily true).

The net result of~(\ref{zzx}) is that an analytic $f$ does not depend upon 
$\bar{z}$, is expandable in a power series in $z$, and this in turn 
justifies the {\em a priori} dubious convention of use of the terminology 
{\em functions of a complex variable}. Indeed, defining the complex
derivative by 
\begin{equation}
\partial_z \equiv \alpha \partial_{x_0} + \beta \partial_{x_1} 
\quad \quad \quad
(~\alpha, \; \beta \in {\cal C}~)~,
\end{equation}
and requiring that
\begin{equation}
\label{cons}
\partial_z z=1 \quad \quad \mbox{and} \quad \quad 
\partial_z \bar{z}=0~,
\end{equation}
we obtain 
\begin{equation}
\label{dz}
\partial_z \equiv \frac{1}{2} \; (\partial_{x_0} - i \partial_{x_1})~.
\end{equation}
The corresponding conjugate complex derivative is 
\[ \partial_{\bar{z}} \equiv \frac{1}{2} \; 
(\partial_{x_0} + i \partial_{x_1})~.
\]
The Cauchy-Riemann condition can then be written in compact form as, 
\begin{equation}
\label{crc}
\partial_{\bar{z}} f(z,\bar{z}) = 0 \quad \quad \quad \left[~
\frac{1}{2} \; (\partial_{x_0} + i \partial_{x_1})~\right]~.
\end{equation}
Thus, these conditions have as a consequence that a mathematical expression 
defining a $z$-differentiable function can depend explicitly only on 
$z=x_0+ix_1$ but not on $\bar{z}=x_0-ix_1$. We also obtain the complex Laplace 
equation 
\begin{equation}
\label{dal}
\Box f(z)=0 \quad \quad \quad (~\Box \equiv \partial_{z} \partial_{\bar{z}}~)~.
\end{equation}

What happens in higher dimensions? Which is the quaternionic
counterpart of $\partial_z$? Are there problems due  
to the noncommutativity of quaternions?

There is a substantial body of literature on attempted constructions
of theories of analytic functions of both real (and even complexified) 
quaternions. In seeking to construct a differential and integral calculus 
of quaternionic functions the first step should be the definition of a 
derivative. A straightforward extension of the complex derivative could be
\[
\partial_{\bar{q}} \equiv 
\frac{1}{2} \left( \partial_{x_0} +
\frac{i\partial_{x_1} + j\partial_{x_2} + k\partial_{x_3}}{3} \right)~,
\]
where 
$\partial_{\bar{q}}$ has been defined in such a way that
\[
\partial_{\bar{q}} q=0 \quad \quad \mbox{and} \quad \quad 
\partial_{\bar{q}} \bar{q}=1~,
\]
in analogy with the above complex case. 
Consequently the Cauchy-Riemann condition would be generalized to 
\begin{equation}
\label{se}
\partial_{x_0} f(q) = 
- \frac{i\partial_{x_1} + j\partial_{x_2} + k\partial_{x_3}}{3} \;  
f(q)~.
\end{equation}
However, due to the noncommutativity of quaternionic numbers, the only
solution to Eq.~(\ref{se}) in the form of a power series of $q$ is
\[ f(q) = c_1 + q c_2 \quad \quad \quad (~c_1, \; c_2 \in {\cal H}~)~.\]
So we find only linear functions in $q$ as solution. Even if these are not 
the only solutions. This is for too 
restrictive to yield any practical use. We also note that~(\ref{se}) does 
{\em not} reduce to~(\ref{zzx}) for complex functions, and indeed does not 
include the standard complex solutions $f(z)$.  

A somewhat different approach extends the concept of holomorphism, which 
for hypercomplex numbers need not coincide with analyticity. 
Working with a noncommutative field we need to define a left or right 
derivative. In fact in the quaternionic version of Eq.~(\ref{derz}) we must 
specify the position of the 
factor $1/\Delta q$. A right quaternionic derivative of the function $f(q)$ 
might be formed by requiring (in analogy to the complex case) that the 
limit
\begin{equation}
\label{derq}
\partial_q  f(q) =  
\lim_{\Delta q \rightarrow 0} ~[f(q+\Delta q) -f(q)]/ \Delta q
\end{equation}
exists and be independent of path for all increments $\Delta q$. By 
considering four linearly independent increments 
$\Delta x_0$, $i\Delta x_1$, $j\Delta x_2$, $k\Delta x_3$ we can derive a 
set of partial differential equations to be satisfied relating the 
components of $f(q)$
\begin{equation}
\label{crqe}
\partial_{x_0} f(q)=
-i\partial_{x_1} f(q)=
-j\partial_{x_2} f(q)=
-k\partial_{x_3} f(q)~.
\end{equation}
This approach also leads to nothing productive since, even for the 
simple function 
$q^2$, the 
\[ \lim_{\Delta q \rightarrow 0} ~\Delta f / \Delta q \]
is not independent of the variation $\Delta q$.  
Even here we do not encompass complex analyticity. Indeed for functions 
independent of both $x_2$ and $x_3$,the above equations allow only
constant functions as solutions.


\section{Fueter Analyticity}
\label{s4}

Among the different approaches to hypercomplex analyticity, the most
important appears to be that of Fueter~\cite{F35} and his school in Z\"urich 
in the 1930's. He 
showed that a third order analytic equation yields 
generalizations of Cauchy's theorem, 
Cauchy's integral formula, and the Laurent expansion~\cite{HD,SUD}. 
Whatismore, the Fueter quaternionic
analyticity~\cite{GUR} has the appreciable virtue of selecting  
something less restrictive than Eqs.~(\ref{se},\ref{crqe}).  

To construct series
with quaternionic coefficients endowed with the ring property we have
to consider the general quaternionic polynomial
\begin{eqnarray}
\label{pol}
p(q) & = & a_0 + a_1 q + q a_2 + a_3 q a_4 + \nonumber \\
     &   & b_1 q^2 + q^2 b_2 + a_3 q^2 a_4 + b_5 q b_6 q b_7 + \nonumber \\\
     &   & . \; . \; .
\end{eqnarray}
which is just a finite sum of quaternionic monomials $m(q)$, e.~g.
\[ m(q) = \alpha_0 q \alpha_2 q ... \alpha_{r-1} q \alpha_r \]
of degree $r$ and with constant quaternions $\alpha_{0,...,r}$. The
polynomial  
$p(q)$ is called general quaternionic since it has no holomorphic
property.  Fueter considered various subgroups of the general 
mapping~(\ref{pol}), the so-called  Weiestrass-like series
\[ L(q) = \sum_n q^n a_n \quad  \mbox{or} \quad R(q) = \sum_n a_n q^n 
\quad \quad \quad [~a_n \in {\cal H}~]~, \]
corresponding, respectively, to left and right holomorphic mappings.

In analogy to the complex case, four derivative operators are introduced 
by
\begin{eqnarray}
{\cal D}_{\bar{q}}^L f(q) & \equiv &
\partial_{x_0} f(q) + i \partial_{x_1} f(q) +
j \partial_{x_2} f(q) +k \partial_{x_3} f(q)~, \nonumber\\ 
{\cal D}_q^L f(q) & \equiv &
\partial_{x_0} f(q) - i \partial_{x_1} f(q) -
j \partial_{x_2} f(q) - k \partial_{x_3} f(q)~,\\
{\cal D}_{\bar{q}}^R f(q) & \equiv &
\partial_{x_0} f(q) + \partial_{x_1} f(q) i +
\partial_{x_2} f(q) j + \partial_{x_3} f(q) k~, \nonumber \\
{\cal D}_q^R f(q) & \equiv &
\partial_{x_0} f(q) - \partial_{x_1} f(q) i -
\partial_{x_2} f(q) j - \partial_{x_3} f(q) k~. 
\end{eqnarray}

The functions
\[ l(q) = {\cal D}_{\bar{q}}^L{\cal D}_q^L L(q)= \Box L(q) \quad \mbox{and} 
\quad r(q) = {\cal D}_{\bar{q}}^R{\cal D}_q^R R(q)= \Box R(q)~,\]
represent, for Fueter, left and right analytic functions since they
are annihilated, respectively by the operators ${\cal D}_{\bar{q}}^L$
and ${\cal D}_{\bar{q}}^R$
\[ {\cal D}_{\bar{q}}^L l(q) = {\cal D}_{\bar{q}}^R r(q) = 0~.\]
So the functions $L(q)$ and $R(q)$ satisfy some Cauchy-Riemann-like
condition, but of the {\em third order} in derivatives instead of the
first. The ``first order'' Laplace equations $\Box L(q)=0$ and $\Box R(q) =0$ do 
not hold in general; however the ``second order'' Laplace equations 
$ \Box^2 L(q)=0$ and $\Box^2 R(q)=0$ are valid. 

Note that given the essential third order nature of the Fueter equations 
there is no simple limit to the complex case. We shall return to this in 
the conclusions.


\section{Global Quaternionic Derivative}
\label{s5}

Returning to first order analytic equations, let us analyze a more general 
case than that of Section III. Consider arbitrary coefficients 
$\in {\cal H}$ of the  ``real'' derivatives 
\begin{equation}
\label{cap}
\partial_q \equiv \alpha \partial_{x_0} + \beta \partial_{x_1} + \gamma
\partial_{x_2} + \delta \partial_{x_3} \quad \quad \quad 
(~\alpha, \; \beta, \gamma, \; \delta \in {\cal H}~)~.
\end{equation}
We find an immediate
problem if we require that our quaternionic derivative satisfies
\begin{equation}
\label{12c}
 \partial_q q = 1 \quad \quad \mbox{and} \quad \quad
\partial_q q^2 = 2q~.
\end{equation}
In fact we obtain the following constraints for the quaternionic
coefficients $\alpha, \; \beta, \; \gamma, \; \delta$
\begin{center}
$\alpha + \beta i + \gamma j + \delta k =1~,$\\
$\alpha i - \beta = i~,~~ 
\alpha j - \gamma = j~,~~ 
\alpha k - \delta = k~.$
\end{center}
Combining the last three equations  we have
\[ 3 \alpha + \beta i + \gamma j + \delta k = 3~,\]
which, when coupled with the first one  gives the trivial
solution
\[ \alpha=1 \]
and consequently
\[ \beta =\gamma  = \delta =0~.\]

Allowing barred quaternionic
coefficients for our quaternionic derivatives, 
a possible solution to the constraints coming from Eq.~(\ref{12c})
is given by 
\begin{center}
$\alpha=1-i\mid i -j\mid j - k \mid k~,$\\
$2\beta = k\mid j - j \mid k~,~~
2\gamma = i\mid k - k \mid i~,~~
2\delta = j\mid i - i \mid j~.$
\end{center}
An immediate check verifies that
\begin{center}
$(\alpha, \; \beta, \; \gamma, \; \delta) \times 1 = 
\left( 4, \; -i, \; -j , \; -k \right)~,$\\
$(\alpha, \; \beta, \; \gamma, \; \delta) \times (i, \; j, \; k) =  
\left( 0, \; -1, \; 0, \; 0 ~;~ 
       0, \; 0, \; -1 , \; 0~;~
       0, \; 0, \; 0, \; -1 \right)~,$
\end{center}
and so the constraints  
on $\alpha, \; \beta, \; \gamma, \; \delta$ are respected. Nevertheless, 
the constraints
coming from the requirement that the quaternionic derivative 
of $q^3$ be $3q^2$ are 
\begin{center}
$\alpha 1+ \beta i + \gamma j + \delta k=1~,$\\
$\alpha i - \beta 1 = i~,~~ 
\alpha j - \gamma 1 = j~,~~
\alpha k - \delta 1 = k~,$\\
$\alpha 1+ \beta i + \frac{1}{3} \; (\gamma j + \delta k) 
=
\alpha 1+ \gamma j  + \frac{1}{3} \; (\beta i + \delta k) 
=
\alpha 1+ \delta k  + \frac{1}{3} \; (\beta i  + \gamma k) 
=1~,$\\
$\beta j + \gamma i =
\beta k + \delta i  = 
\gamma k + \delta j =0~.$
\end{center}
Immediately we obtain 
\begin{center}
$\beta i + \gamma j =
\beta i + \delta k  =
\gamma j + \delta k =0~,$
\end{center}
and so again only the trivial solution exists 
\[ \alpha=1 \quad \quad \mbox{and} 
\quad \quad \beta =\gamma  = \delta =0~.\]

We conclude this Section with the result that it is not possible to obtain 
a satisfactory generalization (e.g. Weiestrass-like series in $q$) 
of the complex
derivative to the quaternionic field {\em if the quaternionic derivative
is defined in terms of \underline{constant} quaternionic (even barred)
coefficients}. In the next Section we will show how the solution arises
by defining ``local'' quaternionic coefficients.


\section{Local Analyticity}
\label{s6}

We rewrite the quaternion $q$ by introducing a ``local'' imaginary 
unit {\em iota}, in the following form
\begin{equation}
q = x_0 + ix_1 + jx_2 + k x_3 \rightarrow x_0 + \iota  x~.
\end{equation}
with
\[
x \equiv |\vec{x}| \quad \quad \quad [~\vec{x}\equiv (x_1, \; x_2, \; x_3)~]
\]
and
\[
\iota  \equiv \frac{ix_1+jx_2+kx_3}{x} ~.
\]
The quaternion $q$ assumes a more familiar form (a complex 
number) and it becomes natural to define the quaternionic derivative in 
terms of the local imaginary unit $\iota$ as
\begin{equation}
\label{local}
\partial_q    =      
\frac{1}{2} \left( \partial_{x_0} -\iota \partial_x \right)~.
\end{equation}
With $\iota$ held fixed (derivative within the $\iota$-complex plane) we 
have, 
\[ 
x\partial_x \equiv 
x_1 \partial_{x_1} + x_2 \partial_{x_2} + x_3 \partial_{x_3}~.
\]
In terms of the global quaternionic imaginary units $i$, $j$ and $k$ 
Eq.~(\ref{local}) reads
\begin{equation}
\partial_q    =      
\frac{1}{2}  \left( \partial_{x_0} - \vec{Q} \cdot \hat{x} \;  
                    \hat{x} \cdot \vec{\partial} \,  
                \right)~,
\end{equation}
where
\begin{eqnarray*}
\vec{Q}         & \equiv & (i, \; j, \; k)~,\\
\vec{\partial}  & \equiv & 
(\partial_{x_1}, \; \partial_{x_2}, \; \partial_{x_3})~,\\
\hat{x}         & \equiv & \vec{x} / x~.
\end{eqnarray*}
This is a ``local'' derivative because the coefficients of 
$\partial / \partial_{1,2,3}$ depend upon the point $q$. To check the 
consistence of our quaternionic derivative it is sufficient to observe
that 
\begin{eqnarray*}
\partial_x (\iota x) & = & \partial_x (ix_1+jx_2+kx_3)\\ 
                       & = & \left(  
     \frac{x_1}{x} \, \partial_{x_1} + 
     \frac{x_2}{x} \, \partial_{x_2} + 
     \frac{x_3}{x} \, \partial_{x_3}  
\right)
      (ix_1+jx_2+kx_3)   \\
                       & = & \frac{ix_1+jx_2+kx_3}{x}\\
                       & = & \iota ~.     
\end{eqnarray*}
We are now able 
 to extend to quaternionic fields the Cauchy-Riemann conditions. 
Explicitly, we find the following quaternionic analyticity requirements:
\begin{equation}
\label{rule}
\partial_{x_0} f(q) = - \iota  \partial_x f(q)~,
\end{equation}
con $f(q)=f_0+if_1+jf_2+kf_3$ ed $f_{0,1,2,3}$ real functions. 
Eq.~(\ref{rule})  can be also rewritten in terms of $(i, \; j, \; k)$ as
follows
\begin{equation}
\partial_{x_0}f_0 + \vec{Q} \cdot \partial_{x_0} \vec{f} =
\hat{x} \cdot {\cal D} \vec{f} 
- \vec{Q} \cdot \left( \hat{x} {\cal D} f_0 + 
\hat{x} \wedge {\cal D} \vec{f} \; \right)~,
\end{equation}
where
\begin{eqnarray*}
\vec{f}         & \equiv & (f_1, \; f_2, \; f_3)~,\\
{\cal D}        & \equiv & \hat{x} \cdot \vec{\partial}~. 
\end{eqnarray*}

The solutions to Eq.~(\ref{rule}) are immediate from complex analysis. 
They consist of all polinomials in $q$ (the ``complex'' $q$) with arbitrary 
right acting quaternionic coefficients. Thus we re-obtain the Fueter 
solutions from a local differential equation. These solutions also contain 
the ``classical'' complex solutions and Eq.~(\ref{rule}) 
reduces to the standard complex equation when applied to a function 
independent of $x_2$ and $x_3$.


\section{Conclusions}

It is useful to recall that one of the difficulties of defining 
quaternionic analyticity is that the simplest extensions to four dimensions 
of complex analyticity tend to restrict the solutions too drastically. 
Another is that some generalizations~(\ref{se},\ref{crqe}) 
do not yield the standard complex limit.

One solution found by Fueter involves a third order analyticity condition 
and yields as solutions a class of polinomials in $q$ with right (or left 
according to the derivative operator used) quaternionic coefficients. 
These solutions $L(q)$ and $R(q)$ are appealing generalizations of the 
standard complex Taylor series, while the associated functions $l(q)$ 
and $r(q)$ satisfy integral properties which are a generalization of 
those of complex functions. This occurs because these latter functions 
satisfy a global first order differential equation, but as a consequence 
they are not simple polinomial series in $q$. Furthermore, while any 
quaternionic function has a complex limit by, say, setting $x_2$ and 
$x_3$ to zero this does not imply that the partial derivatives 
with respect to these missing variables is identically null. Similarly 
while the first order analyticity condition reduce to a complex form, 
the integral equations involving hypersurfaces in $q$-space do not 
automatically become become line integrals in a complex plane. 
Thus $l(q)$ and $r(q)$ do not include complex analytic functions.   

We have shown in this work that the polinomial Fueter solutions $L(q)$ 
and $R(q)$ can be derived from a {\em local} complex analyticity equation, 
i.e. by requiring merely $\iota$-complex analyticity in the natural 
complex plane defined by the point (variable) $q$ and the real axis.
We observe, however, that while $q$ only involves (by definition) 
$\iota$, $f(q)$ is more general because of the unrestricted nature of 
the quaternionic coefficients in its Taylor expansion.

The extensions of our approach to octonions is straightforward 
given the first order nature of our analyticity equation. The corresponding 
definition of $\iota$ is
\[ \iota ~~\equiv ~~\frac{e_1 x_1 + e_2 x_2 + ... + e_7
x_7}{x}~.\]
where $e_1$...$e_7$ are the imaginary unit vectors and $x$ the norm of the 
"vector" part of the octonionic point. The solutions are natural extensions 
of the quaternionic case.


\acknowledgements

We wish to thank D.~Schaulom, R.~Spigler and M.~Vianello for many enlighting 
discussion. We are also grateful to C.~Mariconda for reading the manuscript
and for his useful comments. Finally, for one of us (SdL), it is a pleasure 
to acknowledge F.~De Paolis and P.~Jetzer for their warm hospitality during 
the stay at the Institute for Theoretical Physics, University of Z\"urich, 
where this paper was begun. The work of SdL was supported by Consiglio 
Nazionale delle Ricerche (C.N.R.).



\end{document}